\documentclass[10pt]{amsart}
\usepackage{amssymb}
\usepackage[OT2,OT1]{fontenc}

\newcommand{\ntop}[2]{\genfrac{}{}{0pt}{1}{#1}{#2}}

\let\newpf\proof \let\proof\relax 
\newenvironment{pf}{\newpf[\proofname]}{\qed\endtrivlist}

\def\dim{\operatorname {dim}}

\def\ssigma{\mathfrak{S}_d}
\def\RR{\mathfrak{R}}
\def\RRR {\cal R}

\renewcommand{\P}{\mathbb{P}}

\def\be{\begin{equation}}
\def\ee{\end{equation}}

\def\bm{\begin{pmatrix}}
\def\em{\end{pmatrix}}

\def\ba{{\begin{align}}}
\def\ea{{\end{align}}}

\def\Leb{\mathrm {Leb}}

\def\d{{\underline d}}
\def\l{{\underline l}}

\def\0{{\mathbf 0}}

\def\cal{\mathcal}

\def\GL{\mathrm{GL}}

\newtheorem{thm}{Theorem}[section]
\newtheorem*{thmA}{Theorem A}
\newtheorem*{thmB}{Theorem B}

\newtheorem{lemma}[thm]{Lemma}
\newtheorem{claim}[thm]{Claim}

\theoremstyle{remark}

\numberwithin{equation}{section}

\def \bn {\hfill \\ \smallskip\noindent}

\theoremstyle{definition}

\def\proof{\bn {\bf Proof.} }

\def\note#1
{\marginpar
{\tiny $\leftarrow$
\par
\hfuzz=20pt \hbadness=9000 \hyphenpenalty=-100 \exhyphenpenalty=-100
\pretolerance=-1 \tolerance=9999 \doublehyphendemerits=-100000
\finalhyphendemerits=-100000 \baselineskip=6pt
#1}\hfuzz=1pt}

\newcommand{\dist}{\operatorname{dist}}

\renewcommand{\mod}{\operatorname{mod}}

\newcommand{\CC}{{\cal C}}

\newcommand{\FF}{{\cal F}}

\newcommand{\JJ}{{\cal J}}
\newcommand{\HH}{{\cal H}}

\newcommand{\MM}{{\cal M}}
\newcommand{\NN}{{\cal N}}

\newcommand{\QQ}{{\cal Q}}

\newcommand{\C}{{\mathbb C}}

\newcommand{\N}{{\mathbb N}}

\newcommand{\R}{{\mathbb R}}
\newcommand{\T}{{\mathbb T}}

\newcommand{\Z}{{\mathbb Z}}

\def\B0{{\bold{0}}}

\catcode`\@=12

\def\Empty{}
\newcommand\oplabel[1]{
  \def\OpArg{#1} \ifx \OpArg\Empty {} \else
  	\label{#1}
  \fi}
		
\newcommand{\comm}[1]{}
\newcommand{\comment}[1]{}

\begin{document}

\bigskip\bigskip

\title{Weak mixing for interval exchange
transformations and translation flows}
\author {Artur Avila and Giovanni Forni}

\address{
Laboratoire de Probabilit\'es et Mod\`eles al\'eatoires\\
Universit\'e Pierre et Marie Curie--Boite courrier 188\\
75252--Paris Cedex 05, France
}
\email{artur@ccr.jussieu.fr}
\address{
Topologie et Dynamique\\
Universit\'e Paris-Sud\\
F-91405 Orsay (Cedex), France.
}
\email{Giovanni.Forni@math.u-psud.fr}

\date{\today}

\begin{abstract}
We show that a typical interval exchange transformation is either weakly
mixing or it is an irrational rotation.  We also conclude that a typical
translation flow on a surface of genus $g \geq 2$ (with prescribed
singularity types) is weakly mixing.
\end{abstract}

\setcounter{tocdepth}{1}

\maketitle

\section{Introduction}

Let $d \geq 2$ be a natural number and let $\pi$ be an {\it irreducible} permutation of 
$\{1,\dots,d\}$, 
that is, $\pi \{1,\dots,k\} \neq \{1,\dots,k\}$, $1\leq k<d$.  Given $\lambda\in\R^d_+$, we 
define an 
{\it interval exchange transformation} (i.e.t.) $f:=f(\lambda,\pi)$ in the usual way 
\cite{CFS}, \cite{Ke}: we consider the interval
\be
I:=I(\lambda,\pi)=\left [0,\sum_{i=1}^d \lambda_i \right )\,\,,
\ee
break it into subintervals
\be
I_i:=I_i(\lambda,\pi)=
\left [\sum_{j<i} \lambda_j,\sum_{j \leq i} \lambda_j \right ), \quad
1 \leq i \leq d,
\ee
and rearrange the $I_i$ according to $\pi$ (in the sense that the $i$-th interval 
is mapped onto the $\pi(i)$-th interval). In other words, $f:I \to I$ is given by
\be
x \mapsto x+\sum_{\pi(j)<\pi(i)} \lambda_j - \sum_{j<i} \lambda_j\,, \quad x \in
I_i\,.
\ee
We are interested in the ergodic properties of i.e.t.'s.  Obviously, they
preserve Lebesgue measure. Katok proved that i.e.t.'s and suspension 
flows over i.e.t.'s with roof function of bounded variation are never mixing
\cite{Ka}, \cite{CFS}. Then the  fundamental work of Masur \cite{M} and
Veech \cite{V2} established that almost every i.e.t. is uniquely ergodic
(this means that for every irreducible $\pi$ and for almost every $\lambda 
\in \R^d_+$, $f(\lambda,\pi)$ is uniquely ergodic).

The question of whether the typical i.e.t. is weakly mixing is  more 
delicate except if $\pi$ is a {\it rotation }of $\{1,\dots,d\}$, that is, if $\pi$
satisfies the following conditions: $\pi (i+1) \equiv \pi(i)+1\,\, (\mod\, d)$, for 
all $i\in \{1,\dots,d\}$. In that case  $f(\lambda,\pi)$ is conjugate to a rotation of the circle, 
hence it is not weakly mixing, for every $\lambda \in \R^d_+$.  After the work of Katok 
and Stepin \cite {KS} (who proved weak mixing for almost all i.e.t.'s on  $3$ intervals), 
Veech \cite{V4} established almost  sure weak mixing for infinitely many irreducible 
permutations and asked whether the same property is true for any irreducible permutations  
which is not a rotation.  In this paper, we give an  affirmative answer to this question.

\begin{thmA}

Let $\pi$ be an irreducible permutation of $\{1,\dots ,d\}$ which is not a
rotation.  For almost every $\lambda \in \R^d_+$, $f(\lambda,\pi)$ is weakly
mixing.

\end{thmA}

We should remark that {\it topological} weak mixing has been established
earlier (for almost every i.e.t. which is not a rotation) by Nogueira-Rudolph \cite {NR}.

We recall that a measure preserving transformation $f$ of a probability space $(X,m)$
is said to be {\it weakly mixing } if for every pair of measurable sets $A$, $B\subset X$

\be
\lim_{n\to +\infty}\,\,  \frac{1}{n}\, \sum_{k=0}^{n-1} \vert m(f^{-k} A \cap B) -m(A)m(B)\vert
\,\,=\,\, 0 \,\,.
\ee

It follows immediately from the definitions that  every mixing transformation is weakly mixing
and every weakly mixing transformations is ergodic. A classical theorem states that any 
{\it invertible }measure preserving transformation $f$ is weakly mixing if and only if it has {\it
continuous spectrum}, that is, the only eigenvalue of $f$ is $1$ and the only eigenfunctions
are constants \cite{CFS}, \cite{P}. Thus it is possible to prove weak mixing by ruling out the 
existence of non-constant {\it measurable }eigenfunctions. This is in fact the standard approach
which is also followed in this paper. Topological weak mixing is proved by ruling out the 
existence of non-constant {\it continuous }eigenfunctions. Analogous definitions and 
statements hold for flows.

\subsection{Translation flows} 

Let $M$ be a compact orientable {\it translation surface }of genus $g \geq 1$, that is,
a surface with a finite or empty set $\Sigma$ of conical singularities endowed with an 
atlas such that coordinate changes are given by translation in $\R^{2}$ \cite{GJ1}, 
\cite{GJ2}. Equivalently, $M$ is a compact surface endowed with a flat metric, with 
at most  finitely many conical singularities and trivial holonomy.  For a general flat 
surface  the cone angles at the singularities are $2 \pi (\kappa_1+1) \leq \dots \leq 
2 \pi (\kappa_r+1)$,  where $\kappa_1, \dots ,\kappa_r > -1$ are real numbers 
satisfying $\sum \kappa_i=2g-2$.  If the surface has trivial holonomy, then $\kappa_i
 \in \Z_+$, for all $1 \leq i \leq r$, and there exists a parallel section of the unit tangent 
 bundle $T_{1}M$, that is, a parallel vector field of unit length, well-defined on 
 $M\setminus \Sigma$. A third, equivalent, point of view is to consider pairs $(M,\omega)$ 
 of a compact Riemann surface $M$ and a (non-zero) abelian differential $\omega$. A flat 
 metric on $M$ (with $\Sigma:=\{ \omega=0\}$) is given by $\vert \omega\vert$ and a parallel 
 (horizontal) vector field of unit length is determined by the condition $\omega(v)=1$. 
A {\it translation flow }$F$ on a translation surface $M$ is the flow generated by a parallel 
vector field of unit length on $M\setminus \Sigma$.  The space of all translation flows on a 
given translation surface is naturally identified with the unit tangent space at any regular 
point, hence it is parametrized by the circle $S^{1}$. For all $\theta\in S^{1}$, the
translation flow $F_{\theta}$, generated by the vector field $v_{\theta}$ such that
$e^{-i\theta}\omega (v_{\theta})=1$, coincides with the restriction of the geodesic
flow of the flat metric 
$\vert \omega\vert$ to an invariant surface $M_{\theta}\subset T_{1}M$ (which is the graph 
of the vector field $v_{\theta}$ in the unit tangent bundle over $M\setminus \Sigma$).

\medskip

\noindent The specification of the parameters $\kappa=(\kappa_1,\dots,\kappa_r) \in \Z_+^{r}$
with $\sum \kappa_i=2g-2$ determines a finite dimensional {\it stratum }of the moduli space $\HH(\kappa)$ of translation surfaces which is endowed with a natural complex structure and a Lebesgue measure class \cite{V5}, \cite{Ko}.  We are interested in {\it typical }translation flows (with respect to the Haar measure on $S^{1}$) on {\it typical }translation surfaces (with respect to the Lebesgue measure class on a given stratum). In genus $1$ there are no singularities and translation flows are linear flows on $\T^2$: they are typically uniquely ergodic but never weakly mixing.  In genus $g \geq 2$, the unique ergodicity for a typical translation flow on the typical translation surface was proved by Masur \cite{M} and Veech \cite{V2}.  This result was later strenghtened by Kerckhoff, Masur and Smillie \cite{KMS}  to include arbitrary  translation surfaces. The cohomological equation for a typical translation flow on any translation surface was studied in \cite{F1} where it is proved that, unlike the case of linear flows on the torus, for translation flows on higher genus surfaces there are non-trivial distributional obstructions (which are not measures) to the existence of smooth solutions.  A new proof 
of a similar but finer result on the cohomological equation for typical i.e.t.'s  has been given recently in \cite{MMY1}, \cite{MMY2}. As for interval exchange transformations, the question of weak mixing of translation flows is more delicate than unique ergodicity, but it is widely expected that weak mixing 
holds typically in genus $g \geq 2$. We will show that it is indeed the case:

\comm{
Recall that a measure preserving flow $T_t:M \to M$ is said to be {\it weakly mixing} if there are 
no non-constant measurable solutions $\phi:M \to \R/\Z$ to
\be \label {wm}
\phi(T_t(x))=\phi(x)+ut,
\ee
with $u \in \R$.
}

\begin{thmB}

Let $\HH(\kappa)$ be any stratum of the moduli space of translation surfaces of genus $g \geq 2$.
For almost all translation surfaces $(M,\omega)\in \HH(\kappa)$, the translation flow $F_{\theta}$ on $(M,h)$ is weakly mixing for almost all $\theta\in S^{1}$.

\end{thmB}

\noindent Translation flows and i.e.t.'s are intimately related: the former can be
viewed as suspension flows (of a particular type) over the latter. However,
since the weak mixing property, unlike ergodicity, is not invariant under suspensions
and time changes, the problems of weak mixing for translation flows and
i.e.t.'s are independent of one another. We point out that differently from the case of i.e.t.'s, 
where weak mixing had been proved for infinitely many combinatorics, there had been 
little progress  on weak mixing for typical translation flows (in the measure-theoretic sense), 
except for {\it topological} weak mixing, proved in \cite {L}.  Gutkin and Katok \cite{GK}
proved weak mixing for a $G_{\delta}$-dense set of translation flows on translation surfaces
related to a class of rational polygonal billiards. We should point out that our results tell
us nothing new about the dynamics of rational polygonal billiards (for the well-known 
reason that rational polygonal billiards yield zero measure subsets of the moduli 
space of all translation surfaces).  

\subsection{Parameter exclusion}

To prove our results, we will perform a parameter exclusion to get
rid of undesirable dynamics.  With this in mind, instead or working in 
the direction of understanding the dynamics in the phase space 
(regularity of eigenfunctions\footnote {In this respect, we should remark 
that Yoccoz has pointed out to us the existence
of ``strange'' eigenfunctions for certain values of the parameter.}, etc.),
we will focus on the analysis of the parameter space.

We analyze the parameter space of
suspension flows over i.e.t.'s via a renormalization operator
(the case of i.e.t.'s corresponds to constant roof function). 
This renormalization operator acts non-linearly on i.e.t.'s and linearly on
roof functions, so it has the structure of a cocycle (the Zorich cocycle)
over the renormalization operator on the space of  i.e.t.'s (the Rauzy-Zorich 
induction).
One can work out a criteria for weak mixing (due to Veech \cite{V4})
in terms of the dynamics of the renormalization operator.

An important ingredient in our analysis is the
result of \cite{F2} on the non-uniform hyperbolicity
of the Kontsevich-Zorich cocycle over the Teichm\"uller flow.  This result
is equivalent to hyperbolicity of the Zorich cocycle \cite{Z3}.
(Actually we only need a weaker result, that the Kontsevich-Zorich cocycle, 
or equivalently the Zorich cocycle, has two positive Lyapunov exponents in 
the case of surfaces of genus at least $2$.)

In the case of translation flows a ``linear'' parameter exclusion (on the 
roof function parameters)  shows that ``bad'' roof functions form a small 
set (basically, each positive Lyapunov exponent of the Zorich cocycle gives 
one obstruction for the eigenvalue equation, which has only one free parameter). 
This argument is explained in Appendix A.

The situation for i.e.t.'s is much more complicated, since we have no
freedom of changing the roof function.  We need to do a ``non-linear''
exclusion process, based on a statistical argument. This argument 
proves at once weak mixing  for typical i.e.t.'s and typical translation flows.
While for linear exclusion it is enough to use ergodicity of the
renormalization operator on the space of i.e.t.'s, the statistical argument 
for the non-linear exclusion uses heavily its mixing properties.

\subsection{Outline}

We start this paper with the basic background on cocycles.  We then prove
our key technical result, an abstract parameter exclusion scheme for
``sufficiently random integral cocycles''.

We then discuss known results on the renormalization dynamics for i.e.t.'s and 
show how the problem of weak mixing reduces to the abstract parameter exclusion
theorem.  The same argument also covers the case of translation flows.

In the appendix we present the linear exclusion argument, which is
much simpler than the non-linear exclusion argument but is enough to
deal with translation flows and yields an estimate on the Hausdorff
dimension of the set of translation flows which are not weakly mixing.

\bigskip

{\bf Acknowledgements:} A.A. would like to
thank Jean-Christophe Yoccoz for several very productive discussions
and Jean-Paul Thouvenot for proposing the problem and
for his continuous encouragement. G.F. would like to thank Yakov Sinai
and Jean-Paul Thouvenot who suggested that the results of \cite{F1},
\cite{F2} could be brought to bear on the question of weak mixing for i.e.t.'s.

\section{Background}

\subsection{Strongly expanding maps}

Let $(\Delta,\mu)$ be a probability space. We say that a measurable transformation 
$T:\Delta \to \Delta$, which preserves the measure class of the measure $\mu$, is {\it weakly 
expanding }if there exists a partition (modulo $0$) $\{\Delta^{(l)}, \,l \in \Z\}$ of $\Delta$ 
into sets of positive $\mu$-measure, such that, for all $l \in \Z$, $T$ maps $\Delta^{(l)}$ 
onto $\Delta$, $T^{(l)}:=T \vert \Delta^{(l)}$ is invertible and $T^{(l)}_* \mu$ is equivalent 
to $\mu \vert \Delta^{(l)}$ .

Let $\Omega$ be the set of all finite words with integer entries.  The length (number of entries)
of an element $\l \in \Omega$ will be denoted by $|\l|$.  For any $\l \in \Omega$,
$\l=(l_1,...,l_n)$,
we set $\Delta^\l:=\{x \in \Delta,\, T^{k-1}(x) \in \Delta^{(l_i)},\, 1 \leq k \leq n\}$
and $T^\l:=T^n|\Delta^\l$.  Then $\mu(\Delta^\l)>0$. 

Let $\MM=\{\mu^\l,\, \l \in \Omega\}$, where 
\be
\mu^\l:=\frac {1} {\mu(\Delta^\l)} T^\l_* \mu.
\ee
We say that $T$ is {\it strongly expanding} if
there exists a constant $K>0$ such that
\be
K^{-1} \leq \frac {d\nu} {d\mu} \leq K, \quad \nu \in \MM.
\ee

This has the following consequence.  If $Y \subset \Delta$ is such that
$\mu(Y)>0$ then
\be
K^{-2} \mu(Y) \leq \frac {T^\l_* \nu(Y)} {\mu(\Delta^\l)} \leq K^2 \mu(Y),
\quad \nu \in \MM, \,\l \in \Omega.
\ee

\subsection{Projective transformations}

We let $\P^{p-1}_+ \subset \P^{p-1}$ be the projectivization of $\R^p_+$. 
We will call it the {\it standard simplex}.  A {\it projective contraction}
is a projective transformation taking the standard simplex into
itself.  Thus a projective contraction is the projectivization of some
matrix $B \in \GL(p,\R)$ with non-negative entries.
The image of the standard simplex by a projective contraction is
called a {\it simplex}.  We need the following simple but crucial fact.

\begin{lemma} \label {projectivestronglyexpanding}

Let $\Delta$ be a simplex compactly contained in $\P^{p-1}_+$
and $\{\Delta^{(l)}, \, l \in \Z\}$ be a partition of $\Delta$ (modulo sets 
of Lebesgue measure $0$) into sets of positive Lebesgue measure.  Let 
$T:\Delta \to \Delta$ be a measurable transformation  such that, for all $l \in 
\Z$, $T$ maps $\Delta^{(l)}$ onto $\Delta$, $T^{(l)}:=T \vert \Delta^{(l)}$ is invertible 
 and its inverse is the restriction of a projective contraction.  Then $T$ preserves
a probability measure $\mu$ which is absolutely continuous with respect to Lebesgue 
measure and has a density which is continuous and positive in $\overline \Delta$.  
Moreover, $T$ is strongly expanding with respect to $\mu$.

\end{lemma}

\begin{pf}

Let $d([x],[y])$ be the projective distance between $[x]$ and $[y]$:
\be
d([x],[y])=\sup_{1 \leq i,j \leq p} \left |
\ln \frac {x_i y_j} {x_j y_i} \right |.
\ee
Let $\NN$ be the class of absolutely continuous probability measures on
$\Delta$ whose densities have logarithms which are $1$-Lipschitz with
respect to the projective distance.  Since $\Delta$ has finite projective
diameter, it suffices to show that there exists $\mu \in \NN$ invariant
under $T$ and such that $\mu^\l \in \NN$ for all $\l\in \Omega$.  Notice that $\NN$
is compact in the weak* topology and convex.

Since $(T^\l)^{-1}$ is the projectivization of some matrix $B^{\l}=(b_{ij}^\l)$ in $\GL(p,\R)$ 
with non-negative entries, we have
\be
|\det D(T^\l)^{-1}(x)|=\frac {\|x\|} {\|B^\l \cdot x\|},
\ee
so that
\be
\frac {|\det D(T^\l)^{-1}(y)|} {|\det
D(T^\l)^{-1}(x)|}=\frac {\|B^\l \cdot x\|} {\|B^\l \cdot
y\|} \frac {\|y\|} {\|x\|}=
\left (\frac {\sum_i \left (\sum_j b_{ij}^\l x_j \right )^2}
{\sum_i \left (\sum_j b_{ij}^\l y_j \right )^2} \right )^{1/2} \frac {\|y\|}
{\|x\|} \leq \sup_{1 \leq i \leq p} \frac {x_i \|y\|} {y_i \|x\|} \leq
e^{d([x],[y])}.
\ee
Thus
\be
\Leb^\l:=\frac {1} {\Leb(\Delta^\l)} T^\l_* \Leb \in \NN,
\ee
and
\be
\nu_n:=\frac {1} {n} \sum_{k=0}^{n-1} T_*^k \Leb=\frac {1} {n}
\sum_{\l \in \Omega,|\l|
\leq n} \Leb(\Delta^\l) \Leb^\l \in \NN.
\ee
Let $\mu$ be any limit point of $\{\nu_n\}$ in the weak* topology.  Then $\mu$ is invariant under $T$,
belongs to $\NN$ and, for any $\l \in \Omega$, $\mu^\l$ is a limit of
\be
\nu^\l_n=\left (\sum_{\l^0 \in \Omega,|\l^0| \leq n}
\Leb(\Delta^{\l^0 \l}) \right )^{-1} \sum_{\l^0 \in \Omega,|\l^0|
\leq n} \Leb(\Delta^{\l^0 \l}) \Leb^{\l^0\l} \in \NN,
\ee
which implies that $\mu^\l \in \NN$.
\end{pf}

\subsection{Cocycles} \label {cocycle}

A {\it cocycle} is a pair $(T,A)$ where
$T:\Delta \to \Delta$ and $A:\Delta \to \GL(p,\R)$,
viewed as a linear skew-product $(x,w) \mapsto (T(x),A(x) \cdot w)$ on
$\Delta \times \R^p$.  Notice that $(T,A)^n=(T^n,A_n)$,  where
\be
A_n(x)=A(T^{n-1}(x)) \cdots A(x), \quad n \geq 0.
\ee

If $(\Delta,\mu)$ is a probability space, $\mu$ is an invariant ergodic
measure for $T$ (in particular $T$ is measurable) and
\be
\int_\Delta \ln \|A(x)\| d\mu(x)<\infty,
\ee
we say that $(T,A)$ is a {\it measurable cocycle}.

Let
\be
E^s(x):=\{w \in \R^p,\, \lim \|A_n(x) \cdot w\|=0\},
\ee
\be
E^{cs}(x):=\{w \in \R^p,\, \limsup \|A_n(x) \cdot w\|^{1/n} \leq 1\|.
\ee
Then $E^s(x) \subset E^{cs}(x)$ are subspaces of $\R^p$ (called the {\it
stable} and {\it central stable} spaces respectively), and we have $A(x)
\cdot E^{cs}(x)=E^{cs}(T(x))$, $A(x) \cdot E^s(x)=E^s(T(x))$.  If $(T,A)$ is
a measurable cocycle, $\dim E^s$ and $\dim E^{cs}$
are constant almost everywhere.

If $(T,A)$ is a measurable cocycle, Oseledets Theorem implies that $\lim
\|A_n(x) \cdot w\|^{1/n}$ exists for almost every $x \in \Delta$ and for
every $w \in \R^p$, and that there are
$p$ Lyapunov exponents $\theta_1 \geq \cdots \geq \theta_p$,
characterized by
\begin{multline}
\# \{i,\, \theta_i=\theta\}=
\dim \{w \in \R^p,\, \lim \|A_n(x) \cdot w\|^{1/n} \leq e^\theta\}-\\
-\dim \{w \in \R^p,\, \lim \|A_n(x) \cdot w\|^{1/n}<e^\theta\} \,.
\end{multline}
Thus $\dim E^{cs}(x)=\# \{i,\, \theta_i \leq 0\}$\footnote {It is also possible to 
show that $\dim E^s(x)=\# \{i,\, \theta_i<0\}$.}.
Moreover, if $\lambda<\min \{\theta_i,\, \theta_i>0\}$ then
for almost every $x \in \Delta$, for every subspace $G_0 \subset \R^p$
transverse to $E^{cs}(x)$, there exists $C(x,G_0)>0$ such that
\be \label {g0}
\|A_n(x) \cdot w\| \geq C(x,G_0) e^{\lambda n}\,\| w \|\, , \quad \text{ for all } w \in G_0(x).
\ee

Given $B \in \GL(p,\R)$, we define
\be
\|B\|_0=\max \{\|B\|,\|B^{-1}\|\}.
\ee
If the measurable cocycle $(T,A)$ satisfies the stronger condition
\be
\int_\Delta \ln \|A(x)\|_0 d\mu(x)<\infty
\ee
we will call $(T,A)$ a {\it uniform cocycle}.

\begin{lemma} \label {omega(kappa)}

Let $(T,A)$ be a uniform cocycle and let
\be
\omega(\kappa)=\sup_{\mu(U) \leq \kappa} \sup_{N>0}
\int_U \frac {1} {N} \ln \|A_N(x)\|_0 d\mu(x).
\ee
Then
\be
\lim_{\kappa \to 0} \omega(\kappa)=0.
\ee

\end{lemma}

\begin{pf}

Let
\be \label {omega1}
\omega_N(\kappa)=\sup_{\mu(U) \leq \kappa}
\int_U \frac {1} {N} \ln \|A_N(x)\|_0 d\mu(x).
\ee
Since $\ln \|A_N\|_0$ is integrable for every $N>0$, we have
\be
\lim_{\kappa \to 0} \omega_N(\kappa)=0.
\ee
Let $s=\int_\Delta \ln \|A(x)\|_0 d\mu(x)$.
Let $S_N(x)=\frac {1} {N} \sum_{k=0}^{N-1} \ln \|A(T^k(x))\|_0$, so that
$\frac {1} {N} \ln \|A_N(x)\|_0 \leq S_N(x)$.  We have
\begin{align}
\omega_N(\kappa)&=\sup_{\mu(U) \leq \kappa}
\int_U \frac {1} {N} \ln \|A_N(x)\|_0 d\mu(x) \leq \sup_{\mu(U) \leq \kappa}
\int_U S_N(x) d\mu(x)\\
\nonumber
&\leq \sup_{\mu(U) \leq \kappa}
\int_U s+|S_N(x)-s| d\mu(x)
\leq \kappa s+\|S_N-s\|_{L^1(\mu)}.
\end{align}
By Birkhoff's Theorem, $\|S_N-s\|_{L^1(\mu)} \to 0$, and by (\ref {omega1})
the result follows.
\end{pf}

We say that a cocycle $(T,A)$ is {\it locally constant} if $T:\Delta \to \Delta$
is strongly expanding and $A|\Delta^{(l)}$ is a constant $A^{(l)}$, for all $l\in \Z$.  In this
case,  for all $\l \in \Omega$, $\l=(l_1,...,l_n)$, we let 
\be
A^\l:=A^{(l_n)} \cdots A^{(l_1)}.
\ee

We say that a cocycle $(T,A)$ is {\it integral }if $A(x) \in \GL(p,\Z)$,
for almost all $x \in \Delta$.  An integral cocycle can be regarded as a
skew product on $\Delta \times \R^p/\Z^p$.

\section{Exclusion of the weak-stable space}

Let $(T,A)$ be a cocycle. We define the {\it weak-stable space }at $x\in \Delta$ by
\be
W^s(x)=\{w \in \R^p,\, \|A_n(x) \cdot w\|_{\R^p/\Z^p} \to 0\}\,
\ee
where $\|\cdot \|_{\R^p/\Z^p}$ denotes the euclidean distance from the lattice $\Z^p \subset \R^p$. 
It is immediate to see that, for almost all $x\in \Delta$,  the space $W^s(x)$ is a union of translates of 
$E^s(x)$. If the cocycle is integral, $W^s(x)$ has a natural interpretation as the stable space at $(x,0)$
of the zero section in $\Delta \times \R^p/\Z^p$. If the cocycle is bounded, that is, if the function 
$A:\Delta \to  \GL(p,\R)$ is essentially bounded, then it is easy to see that $W^s(x)=\cup_{c \in
\Z^p} E^s(x)+c$.  In general $W^s(x)$ may be the union of uncountably many translates of $E^s(x)$.

Let $\Theta \subset \P^{p-1}$ be a compact set. We say that $\Theta$ is {\it adapted }to
the cocycle $(T,A)$ if $A^{(l)} \cdot \Theta \subset \Theta$ for all $l$ and if, for
almost every $x \in \Delta$, we have
\be \label {31}
\|A(x) \cdot w\| \geq \|w\|,
\ee
\be \label {32}
\|A_n(x) \cdot w\| \to \infty
\ee
whenever $w \in \R^p \setminus \{0\}$ projectivizes to an element of
$\Theta$.

Let $\JJ=\JJ(\Theta)$ be the set of lines in $\R^p$, parallel to
some element of $\Theta$ and not passing through $0$.

The main result in this section is the following.

\begin{thm} \label {main}

Let $(T,A)$ be a locally constant integral uniform cocycle, and let
$\Theta$ be adapted to $(T,A)$.
Assume that for every line $J \in \JJ:=\JJ(\Theta)$,
$J \cap E^{cs}(x)=\emptyset$ for almost every $x \in \Delta$.
Then if $L$ is a line contained in $\R^p$ parallel to some element of
$\Theta$, $L \cap W^s(x) \subset \Z^p$ for almost every $x \in \Delta$.

\end{thm}

The proof will take the rest of this section.

For $J \in \JJ$, we let $\|J\|$ be the distance between $J$ and $0$.

\begin{lemma}
\label{ss}
There exists $\epsilon_0>0$, such that
\be
\lim_{n \to \infty} \sup_{J \in \JJ}
\mu \left \{x,\,
\ln \frac {\|A_n(x) \cdot J\|} {\|J\|}<\epsilon_0 n \right \} \to 0.
\ee

\end{lemma}

\begin{pf}

Let $C(x,J)$ be the largest real number such that
\be
\|A_n(x) \cdot J\| \geq C(x,J) e^{\lambda n/2} \|J\|, \quad n \geq 0,
\ee
where $\lambda>0$ is smaller than all positive Lyapunov exponents of
$(T,A)$.  By Oseledets Theorem, $C(x,J) \in [0,1]$ is
strictly positive for every $J \in \JJ$
and almost every $x \in \Delta$, and depends continuously on $J$
for almost every $x$.  Thus, 
for every $\delta>0$ and $J \in \JJ$, there exists $C_\delta(J)>0$ such that 
$\mu \{x,\,C(x,J)\le C_\delta(J)\}<\delta$. By Fatou Lemma for any $C>0$ the 
function $F(J):=\mu \{x,\, C(x,J)\le C\}$ is upper semi-continuous, hence 
$\mu \{x,\, C(x,J')\le C_\delta(J)\}<\delta$ for every $J'$ in a neighborhood of 
$J$.  By compactness, there exists $C_\delta>0$ such that
$\mu \{x,\,C(x,J)\le C_\delta\}<\delta$ for every $J \in \JJ$ with
$\|J\|=1$, and hence for every $J \in \JJ$.  The result now follows by
taking $2 \epsilon_0<\lambda$.
\end{pf}

For any $\delta<1/10$, let $W^s_{\delta,n}(x)$ be the set of all $w \in B_\delta(0)$ such that
$\|A_k(x) \cdot w\|_{\R^p/\Z^p}<\delta$ for all $k \leq n$, and let $W^s_\delta(x)=\cap
 W^s_{\delta,n}(x)$.

\begin{lemma}
\label{wss}
There exists $\delta>0$ such that for all $J \in \JJ$ and for almost
every $x \in \Delta$, $J \cap W^s_\delta(x)=\emptyset$.

\end{lemma}

\begin{pf}

For any $\delta<1/10$, let $\phi_\delta(\l,J)$ be the number of connected components of
the set $A^\l(J \cap B_\delta(0)) \cap B_\delta(\Z^p \setminus \{0\})$ and let $\phi_\delta(\l):=
\sup_{J \in \JJ} \phi_\delta(\l,J)$.  For any (fixed) $\l\in \Omega$  the function $\delta \mapsto \phi_\delta(\l)$ is non-decreasing and there exists $\delta_{\l} >0$ such that for $\delta<
\delta_{\l}$ we have $\phi_\delta(\l)=0$. 
We also have
\be \label {35}
\phi_\delta(\l) \leq \|A^\l\|_0 \,, \quad \l \in \Omega\,.
\ee
Given $J$ with $\|J\|<\delta$ and $\l \in \Omega$, let $J_{\l,1},...,J_{\l,\phi_\delta(\l,J)}$ be all the lines
of the form $A^\l \cdot J-c$ where $A^\l(J \cap B_\delta(0)) \cap B_\delta(c) \neq \emptyset$ with $c \in \Z^p \setminus \{0\}$.  Let $J_{\l,0}=A^\l \cdot J$.

By definition we have 
\be
\|J_{\l,k}\|<\delta, \quad k \geq 1.
\ee
To obtain a lower bound we argue as follows. Let $w \in J_{\l,k}$ satisfy $\|w\|=\|J_{\l,k}\|$.  Then $\|w-w'\|<\delta$ for some $w' \in A^\l \,\cdot 
(J \cap B_\delta(0))-c$. Since $J$ is parallel to some element of $\Theta$, it is expanded by 
$A^\l$ (see (\ref {31})). It follows that $\|(A^\l)^{-1} \cdot (w+c)-(A^\l)^{-1} \cdot (w'+c)\|<\delta$, 
which implies $\|(A^\l)^{-1} \cdot (w+c)\|<2 \delta$.  Since $(A^\l)^{-1} \cdot c \in \Z^p
\setminus \{0\}$, we have 
\be
\|A^\l\|_0 \|w\| \geq \|(A^\l)^{-1} \cdot c - (A^\l)^{-1} \cdot (w+c)\|
\geq 1-2\delta,
\ee
and finally we get
\be \label {38}
\|J_{\l,k}\| \geq (1-2\delta) \|A^\l\|^{-1}_0 \geq 2^{-1}\|A^\l\|^{-1}_0\,, \quad k \geq 1.
\ee
On the other hand, it is clear that
\be \label {39}
\|A^\l\|_0 \|J\| \geq \|J_{\l,0}\| \geq \|A^\l\|^{-1}_0 \|J\|.
\ee

Given measurable sets $X,Y \subset \Delta$ such that $\mu(Y)>0$, we let
\be
P_\nu(X|Y)=\frac {\nu(X \cap Y)} {\nu(Y)}, \quad \nu \in \MM,
\ee
\be
P(X|Y)=\sup_{\nu \in \MM} P_\nu(X|Y).
\ee

For $N\in \N\setminus\{0\}$, let $\Omega^N$ be the set of all words of length
$N$, and $\widehat \Omega^N$ be the set of all words of length multiple of $N$.

For any  $0<\eta<1/10$, select a finite set $Z \subset \Omega^N$ such that
$\mu(\cup_{\l \in Z} \Delta^\l)>1-\eta$.  Since the cocycle is locally constant and uniform, there exists
$0<\eta_{0}< 1/10$ such that, for all $\eta<\eta_{0}$, we have
\be \label {312}
\sum_{\l \in \Omega^N \setminus Z} \ln \|A^\l\|_0 \mu(\Delta^\l)<
\frac {1} {2}.
\ee

\begin{claim}
\label{logest}
There exists $N_0\in\N\setminus\{0\}$ such that, if $N>N_0$, then  for every
$Y \subset \Delta$ with $\mu(Y)>0$ we have
\be
\inf_{\nu \in \MM}
\sum_{\l^1 \in Z} \ln \frac {\|J_{\l^1,0}\|} {\|J\|}
P_\nu(\Delta^{\l^1}|\bigcup_{\l \in Z} \Delta^\l \cap T^{-N}(Y)) \geq 2\,.
\ee

\end{claim}

\begin{pf}

By Lemma \ref{ss}, for every $\kappa>0$, there exists $N_0(\kappa)$ such that,
if $N>N_0(\kappa)$, then for every $J \in \JJ$ there exists $Z':=Z'(J) \subset Z$ such that
\be
\ln \frac {\|A^\l \cdot J\|} {\|J\|} \geq \epsilon_0 N,
\ee
\be
\mu\left (\bigcup_{\l \in Z \setminus Z'} \Delta^\l \right )<\kappa.
\ee

We have
\begin{align}
(I)&:=\sum_{\l^1 \in Z'} \ln \frac {\|J_{\l^1,0}\|} {\|J\|}
P_\nu(\Delta^{\l^1}|\bigcup_{\l \in Z} \Delta^\l \cap T^{-N}(Y))\\
\nonumber &\geq
\epsilon_0 N P_\nu(\bigcup_{\l \in Z'} \Delta^\l|\bigcup_{\l \in Z}
\Delta^\l \cap T^{-N}(Y)) \geq
\epsilon_0 N (1-K^4 P_\mu(\bigcup_{\l \in Z \setminus Z'}
\Delta^\l|\bigcup_{\l \in Z} \Delta^\l))\\
\nonumber
&\geq \epsilon_0 N \left (1-K^4 \frac {\kappa} {1-\eta} \right ),
\end{align}
\begin{align}
(II)&:=\sum_{\l^1 \in Z \setminus Z'} \ln \frac {\|J_{\l^1,0}\|} {\|J\|}
P_\nu(\Delta^{\l^1}|\bigcup_{\l \in Z} \Delta^\l \cap T^{-N}(Y))\\
\nonumber
&\geq -\sum_{\l^1 \in Z \setminus Z'} \ln \|A^{\l^1}\|_0
P_\nu(\Delta^{\l^1}|\bigcup_{\l \in Z} \Delta^\l \cap T^{-N}(Y))\\
\nonumber
&\geq -\sum_{\l^1 \in Z \setminus Z'} \ln \|A^{\l^1}\|_0
K^4 P_\mu(\Delta^{\l^1}|\bigcup_{\l \in Z} \Delta^\l)
\geq -K^4 \frac {1} {1-\eta}
\int_{\cup_{\l \in Z \setminus Z'} \Delta^\l} \ln \|A^\l(x)\|
d\mu\\
\nonumber
&\geq -K^4 \frac {1} {1-\eta} \omega(\kappa) N
\end{align}
(where $\omega(\kappa)$ is as in Lemma \ref {omega(kappa)}),
so that for any $\eta<1/10$, for $\kappa>0$ sufficiently small and for all $N>N_0(\kappa)$, we have
\be
\sum_{\l^1 \in Z} \ln \frac {\|J_{\l^1,0}\|} {\|J\|}
P_\nu(\Delta^{\l^1}|\bigcup_{\l \in Z} \Delta^\l \cap T^{-N}(Y)) \geq
(I)+(II) \geq \frac {\epsilon_0} {2} N.
\ee
Hence the claim is proved by taking $N_0 \geq \max \{ N_0(\kappa), 4\epsilon_0^{-1}\}$.
\end{pf}

\begin{claim}

Let $N>N_0$. There exists $\rho_{0}(Z)>0$ such that, for every $0<\rho<\rho_0(Z)$ and every
$Y \subset \Delta$ with $\mu(Y)>0$, we have
\be \label {319}
\sup_{\nu \in \MM}
\sum_{\l^1 \in Z} \|J_{\l^1,0}\|^{-\rho}
P_\nu(\Delta^{\l^1}|\bigcup_{\l \in Z} \Delta^\l \cap T^{-N}(Y)) \leq
(1-\rho) \|J\|^{-\rho}\,.
\ee

\end{claim}

\begin{pf}

Let
\be
\Phi(\nu,Y,\rho):=\sum_{\l^1 \in Z} \frac {\|J_{\l^1,0}\|^{-\rho}}
{\|J\|^{-\rho}}
P_\nu(\Delta^{\l^1}|\bigcup_{\l \in Z} \Delta^\l \cap T^{-N}(Y)).
\ee
Then $\Phi(\nu,Y,0)=1$ and
\be
\frac {d} {d\rho} \Phi(\nu,Y,\rho)=
\sum_{\l^1 \in Z} -\ln \left (\frac {\|J_{\l^1,0}\|}
{\|J\|} \right ) \frac {\|J_{\l^1,0}\|^{-\rho}}
{\|J\|^{-\rho}}
P_\nu(\Delta^{\l^1}|\bigcup_{\l \in Z} \Delta^\l \cap T^{-N}(Y)),
\ee
since $Z$ is a finite set, by Claim  \ref{logest} there exists $\rho_0(Z)>0$ such that,
for every $Y\subset \Delta$ with $\mu(Y)>0$,
\be
\frac {d} {d\rho} \Phi(\nu,Y,\rho)\leq -1, \quad 0 \leq \rho \leq \rho_0(Z),
\ee
which gives the result.
\end{pf}

At this point we fix $0<\eta<\eta_0$, $N>N_0$, $Z\subset \Omega^N$, and $0<\rho<\rho_0(Z)$ so 
that (\ref {312}) and (\ref {319}) hold and let $\delta<1/10$ be so small that we have
\be \label {325}
\sum_{\l \in \Omega^N \setminus Z}
(\rho \ln \|A^\l\|_0+\ln
(1+\|A^\l\|_0 (2 \delta)^\rho)) \mu(\Delta^\l)-
\rho \mu(\bigcup_{\l \in Z} \Delta^\l)=\alpha<0,
\ee
(this is possible by (\ref {312})) 
and
\be \label {324}
\phi_\delta(\l)=0, \quad \l \in Z \,,
\ee
(this is possible since $Z$ is finite).

Let $\Gamma^m_\delta(J)=\{x \in \Delta,\, J \cap W^s_{\delta,m N}(x)
\neq \emptyset\}$.  We must show that $\mu(\Gamma^m_\delta(J)) \to 0$ for every
$J \in \JJ$.  Let
$\psi:\Omega^N \to \Z$ be such that
$\psi(\l)=0$ if $\l \in Z$ and $\psi(\l) \neq \psi(\l')$ whenever $\l \neq
\l'$ and $\l \notin Z$.  We let
$\Psi:\widehat \Omega^N \to \Omega$ be given by
$\Psi(\l^{(1)}\dots \l^{(m)})=\psi(\l^{(1)})\dots \psi(\l^{(m)})$,
where the $\l^{(i)}$ are in $\Omega^N$.
We let $\widehat \Delta^\d=\cup_{\l \in \Psi^{-1}(\d)}
\Delta^\l$.

For $\d \in \Omega$,
let $C(\d) \geq 0$ be the smallest real number such that 
\be 
\sup_{\nu \in \MM} P_{\nu}(\Gamma_{\delta}^m(J)|\widehat \Delta^\d) \leq C(\d) \|J\|^{-\rho},
\quad J \in \JJ.
\ee
It follows that $C(\d) \leq 1$ for all $\d$
(since $\Gamma^m_\delta(J)=\emptyset$, $\|J\|>\delta$).

\begin{claim}

If $\d=(d_1,...,d_m)$, we have
\be \label {341}
C(\d) \leq \prod_{d_i=0} (1-\rho)
\prod_{d_i \neq 0, \psi(\l^i)=d_i} \|A^{\l^i}\|^\rho_0 (1+\|A^{\l^i}\|_0
(2 \delta)^\rho).
\ee

\end{claim}

\begin{pf}

Let $\d=(d_1,...,d_{m+1})$, $\d'=(d_2,...,d_{m+1})$.  There are two
possibilities:
\begin{enumerate}
\item If $d_1=0$, we have by (\ref {319}) and (\ref {324})
\be
P_{\nu}(\Gamma_{\delta}^{m+1}(J)|\widehat \Delta^\d) \leq 
\sum_{\l^1 \in Z}
P(\Gamma_{\delta}^m(J_{\l^1,0})|\widehat \Delta^{\d'})
P_\nu(\Delta^{\l^1}|\widehat \Delta^\d) \leq (1-\rho) C(\d')
\|J\|^{-\rho},
\ee
\item If $d_1 \neq 0$, let
$\l^1$ be given by $\psi(\l^1)=d_1$.  Then either $\|J\|>\delta$ (and
$P(\Gamma^{m+1}_\delta(J)|\widehat
\Delta^\d)=0$) or, by (\ref {35}), (\ref {38}) and (\ref {39}),
\begin{multline}
P(\Gamma_{\delta}^{m+1}(J)|\widehat \Delta^\d) \leq
\sum_{k=0}^{\phi_\delta(\l^1)}
P(\Gamma_{\delta}^m(J_{\l^1,k})|\widehat \Delta^{\d'}) \leq
C(\d') (\|J_{\l^1,0}\|^{-\rho}+\phi_\delta(\l^1) \sup_{k \geq 1}
\|J_{\l^1,k}\|^{-\rho})\\
\leq C(\d') \|J\|^{-\rho} \left (\|A^{\l^1}\|^\rho_0+
\frac {2^\rho \|A^{\l^1}\|^{1+\rho}_0} {\|J\|^{-\rho}} \right )
\leq C(\d') \|J\|^{-\rho} (\|A^{\l^1}\|^\rho_0+             
(2 \delta)^\rho \|A^{\l^1}\|^{1+\rho}_0).
\end{multline}
\end{enumerate}
The result follows.
\end{pf}

Let
\be
\gamma(x):=\left \{ \begin{array}{ll}
-\rho\,, & x \in \cup_{\l \in Z} \Delta^\l\,,\\[5pt]
\rho \ln \|A^\l\|_0+\ln
(1+\|A^\l\|_0 (2 \delta)^\rho), &
x \in \cup_{\l \in \Omega^N \setminus Z} \Delta^\l\,.\end{array}
\right.
\ee
We have chosen $\delta>0$ so that (see (\ref {325}))
\be
\int_\Delta \gamma(x) d\mu(x)=\alpha<0.
\ee

Let $C_m(x)=C(\d)$ for $x \in \widehat \Delta^\d$, $|\d|=m$.  Then by
(\ref {341})
\be
\ln C_m(x) \leq \sum_{k=0}^{m-1} \gamma(T^{k N}(x)) 
\ee
so that, by Birkhoff's ergodic theorem,
$C_m(x) \to 0$ for almost every $x\in\Delta$.  By dominated convergence
(since $C_m(x) \leq 1$), we have
\be
\lim_{m \to \infty} \int_\Delta C_m(x) d\mu(x)=0.
\ee

Notice that
\be
\mu(\Gamma_{\delta}^m(J)) \leq \sum_{\d \in \widehat \Omega, |\d|=m}
\mu(\widehat \Delta^\d)
P_\mu(\Gamma_{\delta}^m(J)|\widehat \Delta^\d) \leq
\int_\Delta C_m(x) \|J\|^{-\rho} d\mu(x),
\ee
so $\lim  \mu(\Gamma_{\delta}^m(J))=0$.
\end{pf}

\noindent{\it Proof of Theorem \ref {main}.}
Assume that there exists a positive measure set $X$ such that for every
$x \in X$, there exists $w(x) \in (L \cap W^s(x)) \setminus \Z^p$.
Thus, for every $\delta>0$ and for every $x \in X$, there exists $n_0(x)>0$ such
that for every $n \geq n_0(x)$, there exists $c_n(x) \in \Z^p\setminus \{0\}$ such that
$A_n(x) \cdot w(x)-c_n(x) \in W^s_\delta(T^n(x))$.

If $A_n(x) \cdot L-c_n(x)$ passes through $0$ for all $n \geq
n_0$, we get a contradiction as follows.
Since $A_n(x)$ expands $L$ (see (\ref {32})) we get   
$\|A_{n-n_0}(T^{n_0}(x))^{-1} (A_n(x) \cdot w(x)-c_n(x))\| \to 0$.
But
\be
A_{n-n_0}(T^{n_0}(x))^{-1}
(A_n(x) \cdot w(x)-c_n(x))=A_{n_0}(x) \cdot w(x)-A_{n-n_0}(T^{n_0}(x))^{-1}
\cdot c_n(x),
\ee
so we conclude that
$A_{n_0}(x) \cdot w(x)=c_{n_0}(x)$, a contradiction.

Thus for every $x \in X$ there exists $n(x) \geq n_0(x)$ such that $A_{n(x)}
\cdot L-c_n(x)$ does not pass through $0$, that is,
$A_{n(x)} \cdot L-c_n(x) \in \JJ$.  By restricting to a subset of
$X$ of positive measure, we may
assume that $n(x)$, $A_{n(x)}(x)$ and $c_{n(x)}(x)$ do not depend on $x\in X$. 
Then $A_{n(x)}(x) \cdot L-c_{n(x)}(x) \in \JJ$ intersects $W^s_\delta(x')$
for all  $x'\in T^{n(x)}(X)$ and $\mu(T^{n(x)}(X))>0$.  This contradicts
Lemma \ref{wss}.
\qed

\section{Renormalization schemes}

Let $d \geq 2$ be a natural number and let $\ssigma$ be the space of {\it irreducible}
permutations on $\{1,...,d\}$, that is $\pi \in \ssigma$ if and only if 
$\pi \{1,...,k\} \neq \{1,...,k\}$ for $1 \leq k < d$.  An i.e.t. $f:=f(\lambda,\pi)$
on $d$ intervals is specified by a pair $(\lambda,\pi) \in \R^d_+ \times \ssigma$ 
as described in the introduction.

\comm{
it is a bijection
of the interval $I=[0,...,\sum \lambda_i)$ obtained by cutting $I$ in
subintervals $I_i=[\sum_{k<i} \lambda_k, \sum_{k \leq i} \lambda_k)$ and
translating them around according to $\pi$: $f(I_i)$ is the $\pi(i)$-th
interval in the new arrangement.
}

\subsection{Rauzy induction}

We recall the definition of the induction procedure first introduced by Rauzy in \cite{R} (see also
Veech \cite{V1}). Let $(\lambda,\pi)$ be such that $\lambda_d \neq \lambda_{\pi^{-1}(d)}$.  Then 
the first return map under $f(\lambda,\pi)$ to the interval
\be
\left [0,\sum_{i=1}^d \lambda_i-\min\{\lambda_{\pi^{-1}(d)},
\lambda_d\} \right )
\ee
can again be seen as an i.e.t. $f(\lambda',\pi')$ on $d$ intervals as follows:
\begin{enumerate}
\item If $\lambda_d<\lambda_{\pi^{-1}(d)}$, let
\be
\lambda'_i=\left \{ \begin{array}{ll}
\lambda_i, & 1 \leq i<\pi^{-1}(d),\\[5pt]
\lambda_{\pi^{-1}(d)}-\lambda_d, & i=\pi^{-1}(d),\\[5pt]
\lambda_d, & i=\pi^{-1}(d)+1,\\[5pt]
\lambda_{i-1}, & \pi^{-1}(d)+1<i \leq d,
\end{array}
\right.
\ee
\be
\pi'(i)=\left \{ \begin{array}{ll}
\pi(i), & 1 \leq i \leq \pi^{-1}(d),\\[5pt]
\pi(d), & i=\pi^{-1}(d)+1,\\[5pt]
\pi(i-1), & \pi^{-1}(d)+1<i \leq d,
\end{array}
\right.
\ee
\item If $\lambda_d>\lambda_{\pi^{-1}(d)}$, let
\be
\lambda'_i=\left \{ \begin{array}{ll}
\lambda_i, & 1 \leq i<d,\\[5pt]
\lambda_d-\lambda_{\pi^{-1}(d)}, & i=d,
\end{array}
\right.
\ee
\be
\pi'(i)=\left \{ \begin{array}{ll}
\pi(i), & 1 \leq \pi(i) \leq \pi(d),\\[5pt]
\pi(i)+1, & \pi(d)<\pi(i)<d,\\[5pt]
\pi(d)+1, & \pi(i)=d,
\end{array}
\right.
\ee
\end{enumerate}
In the first case, we will say that $(\lambda',\pi')$ is obtained from
$(\lambda,\pi)$ by an elementary operation of type $1$, and in the second case
by an elementary operation of type $2$.  In both cases, $\pi'$ is still an irreducible 
permutation.

Let $\QQ_R:\R^d_+ \times \ssigma \to \R^d_+ \times \ssigma$ be the map defined
by $\QQ_R(\lambda,\pi)=(\lambda',\pi')$.  Notice that $\QQ_R$ is defined
almost everywhere (in the complement of finitely many hyperplanes).

The {\it Rauzy class }of a permutation $\pi \in \ssigma$ is the set $\RR(\pi)$ of all $\tilde
\pi$ that can be obtained from $\pi$ by a finite number of elementary operations.
It is a basic fact that the Rauzy classes partition $\ssigma$.

Let $\P^{d-1}_+ \subset \P^{d-1}$ be the projectivization of
$\R^d_+$.  Since $\QQ_R$ commutes with dilations
\be
\QQ_R(\alpha \lambda,\pi)=(\alpha \lambda',\pi')\,, \quad  \alpha \in \R\setminus \{0\}\,,
\ee
$\QQ_R$ projectivizes to a map $\RRR_R:\P^{d-1}_+ \times \ssigma \to \P_+^{d-1}
\times \ssigma$.

\begin{thm}[Masur \cite {M}, Veech \cite{V2}]

Let $\RR\subset \ssigma$ be a Rauzy class. Then $\RRR_R|\P^{d-1}_+
\times \RR$ admits an ergodic conservative infinite absolutely continuous invariant 
measure, unique in its measure class up to a scalar multiple. Its density is 
a positive rational  function. 

\end{thm}

\subsection{Zorich induction}

Zorich \cite{Z1} modified Rauzy induction as follows.  Given $(\lambda,\pi)$, let
$n:=n(\lambda,\pi)$ be such that $Q_R^{n+1}(\lambda,\pi)$ is defined and, for 
$1 \leq i \leq n$, $\QQ_R^i(\lambda,\pi)$ is obtained from $\QQ_R^{i-1}(\lambda,\pi)$ 
by elementary operations of the same type, while
$\QQ_R^{n+1}(\lambda,\pi)$ is obtained from $\QQ_R^n(\lambda,\pi)$ by an
elementary operation of the other type.  Then he sets
\be
\QQ_Z(\lambda,\pi)=\QQ_R^{n(\lambda,\pi)}(\lambda,\pi).
\ee

Then $\QQ_Z:\R^d_+ \times \ssigma \to \R^d_+ \times \ssigma$ is defined almost
everywhere (in the complement of countably many hyperplanes).
We can again consider the projectivization of $\QQ_Z$, denoted by $\RRR_Z$.

\begin{thm}[Zorich \cite {Z1}] \label {invariant measure}

Let $\RR\subset \ssigma$ be a Rauzy class.  Then $\RRR_Z|\P^{d-1}_+
\times \RR$ admits a unique ergodic absolutely continuous probability measure 
$\mu_{\RR}$.  Its density is positive and analytic.

\end{thm}

\subsection{Cocycles}

Let $(\lambda',\pi')$ be obtained from $(\lambda,\pi)$ by the Rauzy or the Zorich induction.  Let $f:=f(\lambda,\pi)$. For any $x\in I':=I(\lambda',\pi')$ and $j\in \{1,\dots, d\}$, let $r_j(x)$ be the 
number of intersections of the  orbit $\{ x,f(x), \dots, f^{k}(x), \dots\}$ with the interval $I_j:= I_j (\lambda,\pi)$ before the first  return time $r(x)$ of $x$ to $I':=I(\lambda',\pi')$, that is, 
$r_j(x):= \# \{0 \le k <r(x) ,  f^{k}(x) \in I_j \}$. In particular, we have $r(x)=\sum_j r_j(x)$.  Notice that $r_j(x)$ is constant on each $I'_i:=I_i(\lambda',\pi')$ and for all $i$, $j\in \{1,\dots, d\}$,  let $r_{ij}:= r_{j}(x)$ for $x\in I'_{i}$. Let $B:=B(\lambda,\pi)$ be the linear operator on $\R^d$ given  by the 
$d\times d$ matrix $(r_{ij})$. The function $B: \P^{d-1}_+ \times \RR \to \GL(d,\R)$ yields a cocycle 
over the Rauzy induction and a related one over  the Zorich induction, called respectively the 
(reduced) {\it Rauzy cocycle }and the (reduced) {\it Zorich cocycle } (denoted respectively by  
$B^R$ and $B^Z$). It is immmediate to see that $B^R,B^Z \in \GL(d,\Z)$, and
\be
B^Z(\lambda,\pi)=B^R(\QQ_R^{n(\lambda,\pi)-1}(\lambda,\pi)) \cdots
B^R(\lambda,\pi).
\ee

Notice that $\QQ(\lambda,\pi)=(\lambda',\pi')$ implies $\lambda=B^*
\lambda'$ ($B^*$ denotes the adjoint of $B$).  Thus
\be \label {lambdalambda'}
\langle \lambda,w \rangle=0 \quad \text {if and only if} \quad
\langle \lambda',B \cdot w \rangle=0.
\ee

Obviously we can projectivize the linear map $B$ and the cocycles $B^R$, $B^Z$.

\begin{thm}[Zorich \cite {Z1}] \label {mu0}

Let $\RR \subset \ssigma$ be a Rauzy class.  We have
\be
\int _{\P^{d-1}_+ \times \RR} \ln \|B^Z\|_0 \,d\mu_{\RR}<\infty \,\,.
\ee

\end{thm}

\subsection{An invariant subspace}

Given a permutation $\pi \in \ssigma$, let  $\sigma$ be the permutation on $\{0,\dots,d\}$ 
defined by
\be
\sigma(i):=\left \{ \begin{array}{ll}
\pi^{-1}(1)-1, & i=0,\\[5pt]
d, & i=\pi^{-1}(d),\\[5pt]
\pi^{-1}(\pi(i)+1), & i \neq 0, \pi^{-1}(d).
\end{array}
\right.
\ee
 For every $j\in \{0,\dots,d\}$, let $S(j)$ be the orbit of $j$ under $\sigma$.  This
defines a partition $\Sigma(\pi):=\{S(j),\, 0 \leq j \leq d\}$ of the set $\{0,...,d\}$.  
For every $S\in \Sigma(\pi)$,  let $b^S \in \R^d$ be the vector defined by
\be
b^S_i:=\chi_S(i-1)-\chi_S(i), \quad 1 \leq i \leq d,
\ee
where $\chi_S$ denotes the characteristic function of $S$.  Let $H(\pi)$ be the annulator of the subspace generated by the set $\Upsilon(\pi):=\{b^S, S \in \Sigma(\pi)\}$. A basic fact from 
\cite{V4} is that if $\QQ(\lambda,\pi):=(\lambda',\pi')$ then
\be
B(\lambda,\pi)^* \cdot \Upsilon(\pi')=\Upsilon(\pi)\,,
\ee
which implies
\be
B(\lambda,\pi) \cdot H(\pi)=H(\pi')\,.
\ee
It follows that the dimension of $H(\pi)$ depends only on the Rauzy class of $\pi\in \ssigma$.  Let 
$N(\pi)$ be the cardinality of the set $\Sigma(\pi)$. Veech showed in \cite{V2} that the
dimension of $H(\pi)$ is equal to $d-N(\pi)+1$ and that  the latter  is  in fact a
non-zero even number equal to $2g$, where $g:=g(\pi)$ is the genus of the Riemann surface obtained by the 
``zippered rectangles'' construction. The space of
zippered rectangles $\Omega(\pi)$ associated to a
permutation $\pi\in \ssigma$ is the space of all
triples $(\lambda, h, a)$ where $\lambda \in  \R^d_+$,
$h$ belongs to a closed convex cone with non-empty
interior $H^+(\pi) \subset H(\pi)$ (specified by finitely many linear inequalities)
and $a$ belongs to a closed parallelepiped $Z(h,\pi) \subset  \overline {\R^d_+}$
of dimension $N(\pi)-1$. Given $\pi \in \ssigma$ and $(\lambda, h, 
a)\in \Omega(\pi)$,  with $h$ in the $H(\pi)$ interior of $H^+(\pi)$,
it  is possible to construct a closed translation surface
$M:=M(\lambda, h, a,\pi)$ of genus $g(\pi)=(d-N(\pi)+1)/2$ by
performing appropriate  gluing operations on the union of the
flat rectangles $R_{i}(\lambda,h) \subset \C$ having bases
$I_{i}(\lambda,\pi)$ and heights $h_{i}$ for $i\in \{1,\dots, d\}$.
The gluing maps are translations specified by the permutation $\pi\in \ssigma$
and by the gluing `heights'  $a:=(a_1,\dots, a_d) \in Z(h,\pi)$. The set
$\Sigma \subset M$ of the singularities of $M$ is in one-to-one correspondence
with the set $\Sigma(\pi)$. In fact for any $S\in \Sigma(\pi)$, the surface $M$
has exactly one conical singularity of total angle $2\pi \nu(S)$, where $\nu(S)$
is the cardinality of $S\cap \{1,\dots, d-1\}$ \cite{V2}. There is a natural local
identification of the relative cohomology $H^{1}(M,\Sigma;\R)$
with the space $\R^d_{+}$ of i.e.t.'s with fixed permutation $\pi\in \ssigma$.
Under this identification  
the generators of $\Upsilon(\pi)$ correspond to integer elements of 
$H^{1}(M,\Sigma;\R)$ and the quotient space of the space generated by $\Upsilon(\pi)$ coincides with the 
absolute cohomology $H^{1}(M,\R)$. By the definition of $H(\pi)$ it follows that $H(\pi)$ is identified with 
the  absolute homology $H_{1}(M,\R)$ and that $H(\pi) \cap \Z^d$ is identified with  $H_1(M,\Z) \subset 
H_1(M, \R)$ (see \cite{Z1}, \S 9), hence $H(\pi) \cap \Z^d$ is a co-compact lattice in $H(\pi)$ and in 
particular we have
\be \label {disth}
\dist(H(\pi),\Z^d \setminus H(\pi))>0.
\ee

\subsection{Lyapunov exponents}

Let $\RR\subset  \ssigma$ be a Rauzy class.  We can consider the restrictions  $B^R(\lambda,\pi)
|H(\pi)$ and $B^Z(\lambda,\pi) |H(\pi)$, $([\lambda],\pi) \in  \P^{d-1}_+ \times \RR$, as integral cocycles 
over $\RRR_Z|\P^{d-1}_+ \times \RR$. We will call these cocycles the Rauzy and Zorich cocycles respectively.  
The Zorich cocycle is uniform (with respect to the measure $\mu_{\RR}$) by Theorems \ref {invariant measure} 
and \ref {mu0}.\footnote {Strictly speaking, to fit into the setting of \S \ref {cocycle} we should fix
an appropriate measurable trivialization of the bundle with fiber $H(\pi)$ at each $(\lambda,\pi)\in 
\P^{d-1}_+ \times \RR(\pi)$ by selecting, for 
each $\tilde \pi \in \RR(\pi)$, an isomorphism $H(\tilde \pi) \to \R^{2g}$ that takes $H(\tilde \pi) \cap 
\Z^d$ to $\Z^{2g}$.}

Let $\theta_1(\RR) \geq \dots  \geq \theta_{2g}(\RR)$ be the Lyapunov exponents of the Zorich cocycle on $\P^{d-1}_+ \times \RR$.  In \cite {Z1}, Zorich showed that $\theta_i(\RR)=-\theta_{2g-i}(\R)$ for all $i\in \{1,\dots, 2g\}$ and that $\theta_1(\RR) >\theta_2(\RR)$ (he derived the latter result from the 
non-uniform 
hyperbolicity of the Teichm\"uller flow proved in \cite{V3}). He also conjectured that $\theta_1(\RR)> \dots 
>\theta_{2g}(\RR)$. Part of this conjecture was proved by the second author in \cite{F2}.

\begin{thm}[Forni, \cite{F2}]

For any Rauzy class $\RR\subset \ssigma$ the Zorich cocycle on 
$\P^{d-1}_+ \times \RR$ is non-uniformly hyperbolic. Thus
\be
\theta_1(\RR)>\theta_2(\RR) \geq \dots \geq
\theta_g(\RR)>0>\theta_{g+1}(\RR) \geq \dots
\geq \theta_{2g-1}(\RR)>\theta_{2g}(\RR).
\ee

\end{thm}

Actually \cite {F2} proved the non-uniform hyperbolicity of a
related cocycle (the Kontsevich-Zorich cocycle), which is a
continuous time version of the Zorich cocycle.  The relation
between the two cocycles can be outlined as follows (see
\cite{V2}, \cite{V3}, \cite{V5}, \cite{Z3}). In \cite{V2} 
Veech introduced a zippered-rectangles ``moduli space''
$\MM(\RR)$ as a quotient of the space  $\Omega(\RR)$ 
of all zippered rectangles associated to permutations
in a given Rauzy class $\RR$. He also introduced a
zippered-rectangles flow on $\Omega(\RR)$ which
projects to a flow on the moduli space $\MM(\RR)$. By 
construction the Rauzy induction is a factor of the return map
of the zippered-rectangles flow to a cross-section
$Y(\RR)\subset \MM(\RR)$.  In fact, such a return map is a
``natural extension'' of the Rauzy induction.  The Rauzy
or Zorich cocycles are cocycles on the bundle with fiber
$H(\pi)$ at $(\lambda,\pi) \in \P^{d-1}\times \RR$.  We recall
that the space $H(\pi)$ can be naturally identified with the
real homology $H_{1}(M,\R)$ of the surface $M:=M(\lambda,h,a,\pi)$.
There is a natural map from the zippered-rectangles
``moduli space'' $\MM(\RR)$ onto a connected component $\CC$ of
a stratum of the moduli space $\HH_{g}$ of holomorphic
(abelian) differentials on Riemann surfaces of genus $g$, and
the zippered-rectangles flow on $\MM(\RR)$ projects onto the
Teichm\"uller flow on $\CC\subset \HH_{g}$.
The Kontsevich-Zorich cocycle, introduced
in \cite{Ko}, is a cocycle over the Teichm\"uller flow on the
{\it real cohomology bundle }over  $\HH_{g}$, that is, the bundle
with fiber the real cohomology $H^{1}(M,\R)$ at every point
$[(M,\omega)] \in \HH_{g}$. The Kontsevich-Zorich cocycle can be 
lifted to a cocycle over the zippered-rectangles flow. The return map 
of the lifted cocycle to  the real cohomology bundle over the
cross-section $Y(\RR)$ projects onto a coycle 
over the Rauzy induction, isomorphic (via Poincar\'e duality)
to the Rauzy cocycle. It follows that the Lyapunov exponents
of the Zorich cocycle on the Rauzy class $\RR$ are related to
the exponents the Kontsevich-Zorich cocycle on $\CC$ 
\cite{Ko}, \cite{F2}, 
\be 
\nu_{1}(\CC)=1 > \nu_{2}(\CC) \geq \dots \nu_{g}(\CC)
>0>\nu_{g+1}(\CC) \geq \dots \geq 
\nu_{2g-1}(\CC)>\nu_{2g}(\CC)=-1.
\ee
by the formula $\nu_{i}(\CC)=\theta_i(\RR)/\theta_1(\RR)$ for all
$i\in \{1,\dots,2g\}$ (see \cite{Z3}, \S4.5).  Thus the non-uniform
hyperbolicity of the Kontsevich-Zorich cocycle on every connected
component of every stratum is equivalent to the hyperbolicity of
the Zorich cocycle on every Rauzy class.

\section{Exclusion of the central stable space for the Zorich cocycle}

\begin{thm} \label {centralstableexclusion}

Let $\pi \in \ssigma$ with $g>1$ and let $L \subset H$ be a line not passing through $0$.  Let 
$E^{cs}$ denote the central stable space of the Rauzy or Zorich cocycle.  If $\dim 
E^{cs}<2g-1$, then for almost every $[\lambda] \in \P^{d-1}_+$, $L \cap E^{cs}([\lambda],\pi)
=\emptyset$.

\end{thm}

Our original proof of this theorem was based on the non-uniform hyperbolicity of the
Kontsevich-Zorich cocycle. The argument is based on the fact that a generic pair of interval 
exchange transformations suspend to a pair $(\cal F, \cal F')$ of {\it transverse }(orientable) 
measured foliations, and their stable spaces (for the Zorich cocycle) can be identified 
with the {\it stable and unstable }spaces (for the Konstevich-Zorich cocycle) at the holomorphic 
differential $\omega\equiv(\cal F, \cal F')$, which are transverse by the non-uniform hyperbolicity.
Yoccoz pointed out to us that it was possible to argue directly with the  Rauzy or Zorich 
cocycle, and later we found out that such an argument was already present in \cite{NR}. The claim 
of \cite{NR} is slightly different from what we need, but the modification is straightforward, so we will 
only give a short sketch of the proof. The advantage of this argument over our original one is that
it does not require the non-uniform hyperbolicity of the cocycle but only $2$ strictly positive
exponents, a property that is much easier to prove than than non-uniform hyperbolicity
(see \cite {F2}).  On the other hand, exclusion of the stable space would use $\dim E^s \leq g$ 
which is an elementary property of the Zorich cocycle.

\begin{pf}

Following Nogueira-Rudolph \cite {NR}, we define $\pi \in \ssigma$ to be standard if
$\pi(1)=n$ and $\pi(n)=1$.  They proved that every Rauzy class contains at least one standard permutation 
\cite {NR}, Lemma 3.2. Clearly it suffices to consider the case when $\pi$ is standard.

Notice that
\be
E^{cs}(\RRR_R([\lambda],\pi))=B^R \cdot E^{cs}([\lambda],\pi)
\ee
for almost every $([\lambda],\pi)$.  It is easy to see (using
Perron-Frobenius together with (\ref {lambdalambda'})) that
$E^{cs}([\lambda],\pi)$ is orthogonal to $\lambda$.

Define vectors $v^{(i)} \in \R^d$ by
\be
v^{(i)}_j=\left \{ \begin{array}{ll}
1, & \pi(j)<\pi(i),\, j>i,\\[5pt]
-1, & \pi(j)>\pi(i),\, j<i,\\[5pt]
0, & \text {otherwise}.
\end{array}
\right.
\ee
It follows that $v^{(i)}$, $1 \leq i \leq d$ generate $H$.

In \cite {NR}, \S 3,
Nogueira and Rudolph showed that  for $1 \leq i \leq
d$ there exist $k_i\in \N$ and a component $D_i \subset \P^{d-1}_+ \times
\{\pi\}$ of the domain of $\RRR_R^{k_i}$ such that $\RRR_R^{k_i}(D_i)=\P^{d-1}_+
\times \{\pi\}$, and defining $B_{(i)}=B^R(\RRR_R^{k_i-1}([\lambda],\pi)) \cdots
B^R([\lambda],\pi)$, we have
\be
B_{(i)} \cdot \bm z_1 \\ \vdots \\ z_d \em=
\bm z_1 \\ \vdots \\ z_d \em+
z_i(v^{(i)}-v^{(d)})-z_d v^{(d)}, \quad i \neq d,
\ee
\be
B_{(d)} \cdot \bm z_1 \\ \vdots \\ z_d \em=
\bm z_1 \\ \vdots \\ z_d \em-z_d v^{(d)}.
\ee

We will now prove the desired statement by contradiction.  If the conclusion of the theorem
does not hold, a density point
argument shows that there exists a set of positive measure of $[\lambda] \in
\P^{d-1}_+$ and a line $L \subset H$ parallel to an element of
$\P^{d-1}_+$ such that
\be
L \cap E^{cs}([\lambda],\pi) \neq \emptyset,
\ee
\be
(B_{(i)} \cdot L) \cap E^{cs}([\lambda],\pi) \neq \emptyset\,, \quad
1 \leq i \leq d.
\ee
Write $L=\{h^{(1)}+th^{(2)},\, t \in \R\}$ with $h^{(1)}$, $h^{(2)} \in H$,
$h^{(2)} \in \R^d_+$
linearly independent.  Then
from $L \cap E^{cs}([\lambda],\pi) \neq \emptyset$, we get
\be
h^{(1)}-\frac {\langle \lambda,h^{(1)} \rangle}
{\langle \lambda,h^{(2)} \rangle}
h^{(2)} \in E^{cs}([\lambda],\pi)\,,
\ee
and similarly we get
\begin{multline}
(h^{(1)}+h_i^{(1)}(v^{(i)}-v^{(d)})-h_{d}^{(1)} v^{(d)})- 
\frac {\langle\lambda,h_1+h_i^{(1)}(v^{(i)}-v^{(d)})-h_d^{(1)} v^{(d)} \rangle }
{\langle\lambda,h^{(2)}+h_i^{(2)}(v^{(i)}-v^{(d)})-h_{d}^{(2)} v^{(d)} \rangle} \times
\\ \times (h^{(2)}+h_i^{(2)}(v^{(i)}-v^{(d)})-h_{d}^{(2)} v^{(d)}) \in E^{cs}([\lambda],\pi)\,,
\end{multline}
for $1 \leq i<d$, and
\be
(h^{(1)}-h_d^{(1)} v^{(d)})-\frac {\langle\lambda,h^{(1)}-h_d^{(1)}
v^{(d)} \rangle} {\langle\lambda, h^{(2)}-h_{d}^{(2)} v^{(d)} \rangle}
(h^{(2)}-h_{d}^{(2)} v^{(d)}) \in E^{cs}([\lambda],\pi)\,.
\ee
A computation then shows that
\be
v^{(1)}, \dots ,v^{(d)} \in E^{cs}([\lambda],\pi)+\{t h^{(2)},\, t \in \R\} \,,
\ee
for almost every such $[\lambda]$.
Thus $E^{cs}([\lambda],\pi)$ has codimension at most $1$ in $H(\pi)$, but this
contradicts $\dim E^{cs}<2g-1=\dim H(\pi)-1$.
\end{pf}

\section{Weak mixing for interval exchange tranformations}

Weak mixing for the interval exchange transformation $f$ is equivalent to the existence of no
non-constant measurable solutions $\phi:I \to \C$ of the equation 
\be
\phi(f(x))=e^{2 \pi i t} \phi(x),
\ee
for any $t \in \R$.  This is equivalent to the following two conditions:
\begin{enumerate}
\item $f$ is ergodic;
\item for any $t \in \R \setminus \Z$, there are no non-zero measurable solutions 
$\phi:I \to \C$ of the equation \be \label {phit} \phi(f(x))=e^{2 \pi i t} \phi(x).
\ee
\end{enumerate}

By \cite {M}, \cite{V2}, the first condition is not an obstruction to almost sure weak mixing: 
$f(\lambda,\pi)$ is ergodic for almost every $\lambda\in \R^d_+$. Our criterion to deal with 
the second condition is the following :

\begin{thm}[Veech, \cite {V4}, \S 7]
\label{veechcriterion}
For any Rauzy class $\RR\subset \ssigma$ there exists an open set $U_{\RR} \subset \P^{d-1}_+ \times \RR$ with the following property.  Assume that the orbit of $([\lambda],\pi) \in \P^{d-1}_+\times \RR$
under  the Rauzy induction $\RRR_R$ visits $U_{\RR} $ infinitely many times. If there exists a non-constant measurable solution $\phi:I \to \C$ to the equation

\be \label {phih} 
\phi(f(x))=e^{2 \pi i t h_i} \phi(x)\, ,  \quad x\in I_{i}(\lambda,\pi)\,,
\ee
with $t \in \R$, $h \in \R^d$, then
\be
\lim_{\ntop{n \to \infty} {\RRR_R^n([\lambda],\pi) \in U_{\RR}}}\,
\|B_n^R([\lambda],\pi) \cdot t h\|_{\R^d/\Z^d}=0.
\ee

\end{thm}

Notice that (\ref {phih}) reduces to (\ref {phit}) when $h=(1,\dots,1)$ and
can thus be used to rule out eigenvalues for i.e.t.'s.  The more general
form (\ref {phih}) will be used in the case of translation flows.

We thank Jean-Christophe Yoccoz for pointing out to us that the above result is due to Veech 
(our original proof does not differ from Veech's). We will call it {\it Veech criterion }for weak mixing.
It has the following consequences:

\begin{thm}[Katok-Stepin, \cite {KS}] \label {katokstepin}

If $g=1$ then either $\pi$ is a rotation or $f(\lambda,\pi)$ is weakly mixing
for almost every $\lambda$.

\end{thm}

(Of course the proof of Katok-Stepin's result precedes Veech criterion.)

\begin{thm}[Veech, \cite{V4}] \label {veech}

Let $\pi\in \ssigma$. If $(1,\dots,1) \notin H(\pi)$, then $f(\lambda,\pi)$ 
is weakly mixing for almost
every $\lambda\in \R^d_+$.

\end{thm}

\subsection{Proof of Theorem A}

By Theorems \ref {katokstepin} and \ref {veech} it is
enough to consider the case where $g>1$ and $(1,...,1) \in H(\pi)$. 
By Veech criterion (Theorem \ref {veechcriterion}), Theorem A is a consequence of the 
following:

\begin{thm} \label {thmA}

Let $\RR\subset \ssigma$ be a Rauzy class with $g>1$, let $\pi \in \RR$ and let $h \in H(\pi) \setminus \{0\}$. Let $U \subset \P^{d-1}_+\times \RR$ be any open set. For almost every $[\lambda] \in \P^{d-1}_+$ 
the following holds: for every $t \in \R$, either $t h \in \Z^d$ or
\be
\limsup_{\ntop {n \to \infty} {\RRR_R^n([\lambda],\pi) \in U}}
\|B_n^R([\lambda],\pi) \cdot t h\|_{\R^d/\Z^d}>0.
\ee

\end{thm}

\begin{pf}

For $n$ sufficiently large there exists a connected component $\Delta \times \{\pi\}  \subset \P^{d-1} \times 
\{\pi\}$ of the domain of $\RRR_Z^n$ which is compactly contained in $U$.  Indeed, the 
connected component of the domain of $\RRR_R^n$ containing $([\lambda],\pi)$ shrinks to 
$([\lambda],\pi)$ as $n \to \infty$, for almost every $[\lambda] \in \P^{d-1}_+$ (this is exactly the 
criterion for unique ergodicity used in \cite{V2}).

If the result does not hold, a density point argument implies that there exists $h \in H(\pi)$ and
a positive measure set of $[\lambda] \in \Delta$ such that
\be \label {criterion1}
\lim_{\ntop {n \to \infty} {\RRR_R^n([\lambda],\pi) \in U}} \|B^R_n
([\lambda],\pi) \cdot th\|_{\R^d/\Z^d}=0, \quad \text {for
some } t \in \R \text { such that } t h \notin \Z^d.
\ee

Let $T:\Delta \to \Delta$ be the map induced by $\RRR_Z$ on $\Delta$.  Then 
$T$ is ergodic, and by Lemma \ref{projectivestronglyexpanding} it is also strongly expanding. For 
almost every $\lambda \in \Delta$, let 
\be 
A(\lambda):=B^Z(T^{r(\lambda)-1}(\lambda),\pi) \cdots B^Z(\lambda,\pi)|H(\pi)\,,
\ee
 where $r(\lambda)$ is the return time of $\lambda \in \Delta$.  Then the cocycle $(T,A)$ is locally 
constant, integral and uniform, and $\Theta:=\overline {\P^{d-1}_+}$ is adapted to $(T,A)$. The central 
stable space of $(T,A)$ coincides with the central stable space of $(\RRR_Z, B^Z|H(\pi))$ almost everywhere.  
Using Theorem \ref{centralstableexclusion}, we see that all the hypotheses of Theorem \ref {main} are 
satisfied.  Thus for almost every $[\lambda] \in \Delta$, the line $L=\{t h,\, t \in \R\}$ intersects the 
weak stable space in a subset of $H(\pi) \cap \Z^d$.  This implies (together with (\ref {disth})) that (\ref 
{criterion1}) fails for almost every $\lambda\in \Delta$, as required.
\end{pf}

\section{Translation flows}

\subsection{Special flows}

Any translation flow on a translation surface can be regarded, by considering its return map to a transverse interval, as a {\it special flow} (suspension flow) over some interval exchange 
transformation with a roof function constant on each sub-interval. For completeness we discuss 
weak mixing for general special flows over i.e.t.'s with sufficiently regular roof function. Thanks 
to recent results on the cohomological equation  for i.e.t.'s \cite{MMY1}, \cite{MMY2} the general 
case can be reduced to the case of special flows with roof function constant on each sub-interval. 

Let $F:=F(\lambda,h,\pi)$ be the special flow over the i.e.t.  $f:=f(\lambda,\pi)$ with roof function 
specified by the vector $h\in \R^d_+$, that is, the roof function is constant equal to $h_{i}$ on the
sub-interval $I_{i}:= I_{i}(\lambda,\pi)$, for all $i\in \{1,\dots,d \}$. We remark that, by Veech's 
``zippered rectangles'' construction (see \S 6), if $F$ is a translation flow then necessarily $h\in H(\pi)$. 

The phase space of $F$  is the union of disjoint rectangles $D:=\cup_i I_i \times [0,h_i)$, and the flow $F$ 
is completely determined by the conditions $F_s(x,0)=(x,s)$, $x \in I_i$, $s<h_i$, $F_{h_i}(x,0) =(f(x),0)$, 
for all $i\in \{1,\dots,d\}$. Weak mixing for the flow $F$ is equivalent to the existence of no 
non-constant measurable solutions $\phi:D \to \C$ of the equation
\be
\phi(F_s(x))=e^{2 \pi i t s} \phi(x),
\ee
for any $t \in \R$, or,  in terms of the i.e.t. $f$, 
\begin{enumerate}
\item $f$ is ergodic;
\item  for any $t \neq 0$ there are no non-zero measurable solutions $\phi:I \to \C$ of
equation (\ref {phih}).
\end{enumerate}

\begin{thm} \label {thmB}

Let $\pi \in \ssigma$ with $g>1$.
For almost every $(\lambda,h) \in \R^d_+ \times (H(\pi) \cap \R^d_+)$,
the special flow $F:=F(\lambda,h,\pi)$ is weakly mixing.

\end{thm}

\begin{pf}

This is an immediate consequence of Veech criterion and Theorem \ref{thmA}.
\end{pf}

This theorem is all we need in the case of translation flows since it takes care of the case 
$h\in H(\pi)$. Let $H^{\perp}(\pi)$ be the orthogonal complement of $H(\pi)$ in $\R^{d}_{+}$. 
The case when $h\in \R^{d}_{+}$ has non-zero orthogonal projection on $H^{\perp}(\pi)$
is covered by the following:

\begin{lemma}[Veech, \cite {V3}]
\label{Veechlemma}
Assume that $([\lambda],\pi) \in \P^d_+ \times \ssigma$ is such that
$\QQ_R^n([\lambda],\pi)$ is defined for all $n>0$.  If for some $h \in \R^d$
and $t \in \R$, we have 
\be
\liminf_{n \to \infty}
\|B_n^R([\lambda],\pi) \cdot t h\|_{\R^d/\Z^d}=0\,,
\ee
then the orthogonal projection of $t h$ on
$H^{\perp}(\pi)$ belongs to $\Z^d$.

\end{lemma}

This lemma, together with Veech criterion, can be used to establish
typical weak mixing for special flows with some specific combinatorics
(such that in particular $\dim H(\pi) \leq d-2$). However, it does not help 
at all when  $h \in H(\pi)$, which is the relevant case for translation flows.

\begin{thm}
\label{wm}
Let $\pi \in \ssigma$ with $g>1$ and let $h \in \R^d \setminus \{0\}$.
If $U \subset \P^{d-1}_+$ is any open set, then for almost every $[\lambda] \in \P^{d-1}_+$ and 
for every $t \in \R$, either $t h \in \Z^d$ or
\be \label {th}
\limsup_{\RRR_R^n([\lambda],\pi) \in U}
\|B_{n}^R([\lambda],\pi) \cdot t h\|_{\R^d/\Z^d}>0.
\ee

\end{thm}

\begin{pf}

If $h \in H^{\perp}(\pi)$, this is just a consequence of Lemma \ref{Veechlemma}.  So we
assume that $h \notin H^{\perp}(\pi)$. Let $w$ be the orthogonal projection of $h$ on $H(\pi)$.  By Theorem \ref {thmA},  there exists a full measure set of $[\lambda] \in \P^{d-1}_+$ such that if 
$t w \notin \Z^d$
then
\be
\limsup_{\RRR_R^n([\lambda],\pi) \in U}
\|B_{n}^R([\lambda],\pi) \cdot t w\|_{\R^d/\Z^d}>0.
\ee
By Lemma \ref{Veechlemma}, if (\ref {th}) does not hold for some $t\in \R$, then $t
h=c+tw$ with $c \in \Z^d$.  But this implies that
\be
\|B_{n}^R([\lambda],\pi) \cdot t w\|_{\R^d/\Z^d}=
\|B_{n}^R([\lambda],\pi) \cdot t h\|_{\R^d/\Z^d},
\ee
and the result follows.
\end{pf}

\begin{thm} \label {cor}
Let $\pi \in \ssigma$ with $g>1$.
For almost every $(\lambda,h) \in \R^d_+ \times \R^d_+$,
the special flow $F:=F(\lambda,h,\pi)$ is weakly mixing.

\end{thm}

\begin{pf}

It follows immediately from Theorem \ref{wm} and Veech
criterion (Theorem \ref{veechcriterion}) by Fubini's theorem.
\end{pf}

Following  \cite{MMY1},  \cite{MMY2}, we let $BV(\sqcup I_{i})$ be the space of functions whose restrictions to each of the intervals $I_{i}$  is a function of bounded variation,  $BV_{\ast}(\sqcup I_{i})$ be the hyperplane of $BV(\sqcup I_{i})$ made of functions whose integral on the disjoint union
$\sqcup I_{i}$ vanishes and $BV_{\ast}^{1}(\sqcup I_{i})$ be the space of 
functions which are
absolutely continuous on each $I_{i}$ and whose first derivative belongs 
to $BV_{\ast}
(\sqcup I_{i})$. Let $\chi:BV_{\ast}^{1}(\sqcup I_{i})\to \R^{d}$ be the
surjective linear map defined by
\be
\chi(r):= \left(\int_{I_1} r(x)dx\, , \dots, \int_{I_d} r(x)dx\right)\,,
\quad r\in BV_{\ast}^{1}(\sqcup I_{i})\,.
\ee

\begin{thm}

Let $\pi \in \ssigma$ with $g>1$. For almost every $\lambda \in \R^d_+$, there 
exists a full measure set $\FF\subset \R^{d}$ such that if
$r\in BV_{\ast}^{1}(\sqcup I_{i})$ is a strictly positive function 
with $\chi(r)\in \FF$, 
then the special  flow $F:=F(\lambda,\pi;r)$ over the i.e.t. 
$f:=f(\lambda,\pi)$ under the  roof function 
$r$ is weakly mixing.
 
 \end{thm}

\begin{pf} 
By the definition of a special flow over the map $f$ and under the roof function $r$ (see \cite{CFS},
Chap. 11),  the flow $F$ has continuous spectrum if and only if 

\begin{enumerate}
 \item $f$ is ergodic, 
 \item  for any $t\neq 0$ there are no non-zero measurable solutions $\phi:I \to \C$ of the equation
 \be
 \label{phir}
 \phi(f(x)) = e^{2\pi i t r(x)} \phi(x)\,, \quad x\in I\,.
 \ee
 \end{enumerate}
By \cite{MMY1}, \cite{MMY2}, under a full measure condition on $\lambda 
\in \R^d_+$ (a Roth-type condition) the cohomological equation
\be
u(f(x))-u(x) = r(x) - \chi_{i}(r)\,,  \quad x\in I_{i}\,,
\ee
has a bounded measurable solution $u:I \to \R$ (the Roth-type condition 
also implies that $f$ is
uniquely ergodic). If for some $t\neq 0$ there exists a solution $\phi$ of equation (\ref{phir}), then the function $\psi:I\to \C$ given by 
\be
\psi(x):=e^{-2\pi i t  u(x)} \phi(x)\,, \quad x\in I\,,
\ee
is a solution of the equation (\ref{phih}) for $h=\chi(r)\in \R^{d}$.
The result then follows from Theorem \ref{wm} and from Veech 
criterion (Theorem \ref{veechcriterion}) by Fubini's theorem.
\end{pf}
 
\subsection{Proof of Theorem B}

Let $\CC$ be a connected component of a stratum of the moduli space
of holomorphic differentials of genus $g>1$. By
Veech's ``zippered rectangles'' construction $\CC$ can be 
locally parametrized by triples $(\lambda,h,a) \in \Omega(\pi)$
where $\pi$ is some irreducible permutation with $g(\pi)=g$
(see \cite{V2}, \cite{V3}). Moreover, this parametrization 
(which preserves the Lebesgue measure class) is such that the special 
flow $F:=F(\lambda,h,\pi)$ 
is isomorphic to the vertical translation flow on the translation surface $M(\lambda,h,a,\pi)$.  Thus Theorem B follows from Theorem \ref {thmB} by Fubini's theorem.

\appendix
\section{Linear exclusion}

\begin{thm}
\label{HD}
Let $(T,A)$ be a measurable cocycle on $\Delta\times \R^{p}$.  For almost every $x \in \Delta$, if $G \subset \R^p$ is any affine subspace parallel to a linear subspace $G_{0}\subset \R^p$ transverse to 
the central stable space $E^{cs}(x)$, then  the Hausdorff dimension of $G \cap W^s(x)$ is equal to $0$.

\end{thm}

\begin{pf}

Let $x\in \Delta$. If $n \geq m \geq 0$ we let $S_{\delta,m,n}(x)$ be the set of
$w \in G$ such that $A_k(x) \cdot w \in B_\delta(\Z^p)$, for all $m \leq k \leq n$. 
Thus
\be
G \cap W^s(x)=\cap_{\delta>0} \cup_{m \geq 0} \cap_{n \geq m}
S_{\delta,m,n}(x).
\ee
If $\delta<1/2$, all connected components of $S_{\delta,m,n}(x)$ are convex open sets of diameter at most
\be
2 \delta \,C(x,G_0)^{-1} e^{-\lambda n},
\ee
where $C(x,G_0)>0$ and $\lambda>0$ are given by Oseledets Theorem as in (\ref {g0}).
For $n \geq 0$, let $\rho_n(x)$ be the maximal number of connected components of 
$S_{\delta,m,n+1}(x)$ intersecting $U$, over all $m \leq n$ and all connected components 
$U$ of $S_{\delta,m,n}(x)$.  We have 
\be
\rho_n(x) \leq 1+(3\delta \|A(T^n(x))\|)^p.
\ee
Let then
\be
\beta_{\delta}(x):=
\limsup_{n \to \infty}
\frac {1} {n} \sum _{k=0}^{n-1}\ln (1+(3 \delta \|A(T^k(x))\|)^p).
\ee
By Birkhoff's ergodic theorem, $\lim_{\delta \to 0}\beta_{\delta}(x)=0$ for almost every $x \in \Delta$.

It follows that there exists a sequence $\epsilon_{\delta}(x,n)$, with $\lim_{n \to \infty}
\epsilon_{\delta}(x,n)=0$ for almost every $x\in \Delta$, such that each connected component $U$ of $S_{\delta,m,m}(x)$ intersects at most
\be
\prod_{k=m}^{n-1} \rho_k(x)  \leq 
e^{\epsilon_{\delta}(x,n) n} \,e^{\beta_{\delta}(x) n}
\ee
connected components of $S_{\delta,m,n}(x)$. Thus $U$ intersects $\cap_{n \geq m} S_{\delta,m,n}(x)$ in a set of upper box dimension at most $\frac {\beta_{\delta}(x)} {\lambda}$.  We
conclude that $\cup_{m \geq 0} \cap_{n \geq m} S_{\delta,m,n}(x)$
has Hausdorff dimension at most $\frac {\beta_{\delta}(x)} {\lambda}$, hence $G \cap W^s(x)$ 
has Hausdorff dimension $0$, for almost every $x \in \Delta$.
\end{pf}

Theorem \ref{HD}, together with Veech's criterion (Theorem \ref{veechcriterion}) has the following consequence, which implies Theorem B.

\begin{thm}

Let $\pi \in \ssigma$.  Then for almost every $\lambda \in
\R^d_+$, the set of $h \in H(\pi) \cap \R^d_+$ such that $(\lambda,h,\pi)$ is not
weakly mixing has Hausdorff dimension at most $g(\pi)+1$.\footnote {In order to get
a weaker, measure zero, statement (when $g(\pi)>1$), it is enough to assume as above  
that the Zorich cocycle has two positive Lyapunov exponents.}

\end{thm}

\begin{pf}

It is enough to show that, for almost every $x\in \Delta$, the weak stable space $W^s(x)$ of the cocycle $(T,A)$ considered in Theorem \ref {thmA} has Hausdorff dimension at most $g(\pi)$. In fact, in this
case the set of $h \in H(\pi) \setminus \{0\}$ such that the line through $h$ intersects $W^s(x)$ in some $w \neq 0$ has Hausdorff dimension at most $g(\pi)+1$, and the result then follows by Veech's criterion.  Since $W^s(x)=(G_0 \cap W^s(x))+E^s(x)$, where $G_{0}$ is any linear subspace transverse to the stable space $E^s(x)$, it is enough to show that $W^s(x) \cap G_0$ has Hausdorff dimension $0$. This follows from the above Theorem \ref{HD}, since by  the non-uniform hyperbolicity of $(T,A)$ the central stable space $E^{cs}(x)$ and the stable space $E^{s}(x)$ coincide, for  almost every $x\in \Delta$.
\end{pf}

\end{document}